\documentclass[11pt]{amsart}

\usepackage{epsf,amssymb}

\newtheorem{thm}{Theorem}[section]
\newtheorem{prop}[thm]{Proposition}

\newtheorem{lem}[thm]{Lemma}

\newtheorem{question}{Question}

\newcommand{\R}{\mbox{\bf{R}}}
\newcommand{\Z}{\mbox{\bf{Z}}}
\newcommand{\Q}{\mbox{\bf{Q}}}

\newcommand{\bdry}{\partial}

\newcommand{\tb}{tb}
\newcommand{\tw}{t}
\newcommand{\s}{\vskip.1in}
\newcommand{\n}{\noindent}

\newcommand{\be}{\begin{enumerate}}
\newcommand{\ee}{\end{enumerate}}

\begin{document}
\title{Tight contact structures with no symplectic fillings}

\author{John B. Etnyre}
\address{Stanford University, Stanford, CA 94305}
\email{etnyre@math.stanford.edu}
\urladdr{http://math.stanford.edu/\char126 etnyre}

\author{Ko Honda}
\address{University of Georgia, Athens, GA 30602}
\email{honda@math.uga.edu}
\urladdr{http://www.math.uga.edu/\char126 honda}

\thanks{JE supported in part by NSF Grant \# DMS-9705949.  KH supported in part by
NSF Grant \# DMS-0072853 and the American Institute of Mathematics.}

\date{First version: October 4, 2000.  This version: October 17, 2000.}


\begin{abstract}
	We exhibit tight contact structures on 3-manifolds that
	do not admit any symplectic fillings.
\end{abstract}
\maketitle

\section{Introduction}


In the early 1980's, D. Bennequin \cite{Bennequin} proved the existence of {\it exotic} contact
structures on $\R^3$.  These were obtained from the {\it standard} contact structure on $\R^3$ given by
the 1-form $\alpha=dz-ydx$, by performing modifications called {\it Lutz twists}.   The key
distinguishing feature was that the exotic contact structures contained {\it overtwisted
disks}, i.e., disks $D$ which are everywhere tangent to the 2-plane field distribution along
$\bdry D$.    On the
other hand, using an ingenious argument which used braid foliations, Bennequin succeeded in proving that
the standard contact structure contained no overtwisted disks.  In a strange twist of fate, the exotic
contact structures eventually turned out not to be so exotic, when Eliashberg \cite{Eliashberg89} gave a
complete classification of contact structures which contain overtwisted disks (called {\it overtwisted}
contact structures) in terms of homotopy theory.

With the advent of Gromov's
theory of holomorphic curves \cite{Gromov}, it became easier to determine when a contact structure on a
3-manifold is {\it tight}, i.e., contains no overtwisted disks \cite{Eliashberg90a}.   Loosely
speaking, a contact structure is {\it symplectically fillable} if it is the boundary of a
symplectic 4-manifold. Gromov and Eliashberg showed that a symplectically fillable structure
was necessarily tight.  In fact, until the mid-1990's, all known tight contact structures
were shown to be tight using symplectic fillings. This included two rich sources of tight
contact structures --- perturbations of taut foliations as in \cite{et} and Legendrian surgery
as in \cite{Eliashberg90} and \cite{Weinstein91}.   
This promoted 
Eliashberg and others to ask whether tight contact structures are the same as symplectically 
fillable contact structures. Subsequently, gluing techniques were developed by Colin 
\cite{Co97, Co99} and Makar-Limanov \cite{ML}, and strengthened in \cite{Honda4}.   Largely due 
to the improvements in gluing techniques, tight contact structures could now be constructed 
without resorting to symplectic filling techniques. The main result of this paper shows that 
the symplectically fillable contact structures form a proper subset of tight contact 
structures.

\begin{thm}\label{main}
	Let $M_1$ (resp. $M_2$) be the Seifert fibered space over $S^2$ with Seifert invariants $(-\frac{1}{2},
	\frac{1}{4},\frac{1}{4})$ (resp. $(-\frac{2}{3}, \frac{1}{3},\frac{1}{3})$).  Then $M_1$ admits one tight contact structure
	and $M_2$ admits two nonisotopic tight contact structures that are not weakly symplectically semi-fillable.
\end{thm} 

In this paper we will provide a complete proof for $M=M_1$; the proof for $M=M_2$ is similar, and 
we will briefly discuss the necessary modifications at the end of Section \ref{mainresult}.

\s\n
{\it Remark on notation.}  A Seifert fibered space over a closed oriented surface $\Sigma$
with $n$ singular fibers is often denoted by
$(g; (1,e), (\alpha_1,\beta_1),\cdots, (\alpha_n,\beta_n))$, or by
$(g; e, {\beta_1\over \alpha_1},$ $\cdots, {\beta_n\over \alpha_n})$, where $g$ is the genus of the base $\Sigma$,
$e\in \Z$ is the Euler number,
and $\alpha_i,\beta_i\in \Z^+$ are relatively prime.   In this notation, $(-{1\over 2}, {1\over 4}, {1\over 4})$
would correspond to $(0; -1, {1\over 2}, {1\over 4}, {1\over 4})$ and
$(-{2\over 3}, {1\over 3}, {1\over 3})$ to $(0; -1, {1\over 3}, {1\over 3}, {1\over 3})$.

\section{Background and preliminary notions}

We briefly review the basic notions in Section \ref{first} and proceed to a
discussion of symplectic fillings in Section~\ref{fillings}.  There we introduce the various types of
symplectic fillings and discuss the work of Lisca concerning the non-existence of
fillable structures on certain manifolds. Finally, in Section \ref{leg}, we discuss the contact surgery
technique which is usually called Legendrian surgery.

Convex surface theory will be our main tool throughout this paper. Originally developed by
Giroux in \cite{Giroux91}, there have been many recent papers discussing convex surfaces. All
the facts relevant to this paper concerning convex surfaces may be found in \cite{Honda1, EH}
(see also \cite{Kanda, Honda2}), and we assume the reader is familiar with the terminology in
these papers.

\subsection{Contact structures and Legendrian knots}  \label{first}
In this section we review a few basic notions of contact topology in dimension three.
This is more to establish terminology than to introduce the readers to these
ideas. Readers unfamiliar with these ideas should see \cite{a, Eliashberg92}.

An oriented 2-plane field distribution $\xi$ on an oriented 3-manifold $M$ is called a
{\it positive contact structure} if $\xi=\ker \alpha$ for some global 1-form $\alpha$
satisfying $\alpha\wedge d\alpha >0$.  The 1-form $\alpha$ is called the {\it contact form}
for $\xi$.  In this paper we will always assume that our ambient manifold $M$ is oriented and
the contact structure $\xi$ is positive and oriented.
If $\Sigma$ is a surface in a contact manifold $(M,\xi)$, then $\Sigma$ has a singular
foliation $\Sigma_\xi,$ called the {\em characteristic foliation}, given by integrating the
singular line filed $T_x\Sigma\cap \xi_x.$
It is important to remember that the characteristic foliation on a surface determines the germ of the contact
structure along the surface.
A contact structure $\xi$ is called {\em overtwisted} if there is an embedded disk $D$ which is
everywhere tangent to $\xi$ along $\bdry D$. A contact structure is
{\em tight} if it is not overtwisted. For a complete classification of overtwisted contact
structures see \cite{Eliashberg89}. A tight contact structure $\xi$ on $M$ is called {\em
universally tight} if  it remains tight when pulled back to the universal cover, and is called
{\em virtually overtwisted} if it becomes overtwisted when pulled back to some finite cover. It is an
interesting problem to determine whether every tight contact structure is either universally
tight or virtually overtwisted.

A knot $K$ embedded in a contact
manifold $(M,\xi)$ is called {\em Legendrian} if it is everywhere tangent to $\xi.$   A choice of
nonzero section of $\xi$ transverse to $K$ gives a framing of the normal bundle of $K$, usually called
the {\it contact framing}. If $\mathcal{F}$ is some preassigned framing of $K$, then we
associate an integer called the {\em twisting number of $\xi$ along $K$} relative to
$\mathcal{F}$, which is the difference in twisting between the contact framing and $\mathcal{F}$,  and
denote it $\tw(K,\mathcal{F})$ (or just $\tw(K)$ if the framing $\mathcal{F}$ is
understood). If $K$ is a null-homologous knot and $\mathcal{F}$ is given by a Seifert surface for $K$,
then $\tw(K)$ is called the {\em Thurston-Bennequin invariant} of $K$ and is usually denoted
$\tb(K).$

A closed surface or a properly embedded compact surface $\Sigma$ 
with Legendrian boundary is called {\it convex} if there exists
a contact vector field everywhere transverse to $\Sigma$.  To a convex surface 
$\Sigma$ we associate an isotopy class of multicurves called the 
{\it dividing set} $\Gamma_\Sigma$ (or simply $\Gamma$).  If $\Sigma$ is closed, then components of
$\Gamma_\Sigma$ are closed curves, and if $\Sigma$ has boundary, there may also be properly embedded
arcs.  The number of components of $\Gamma_\Sigma$ is written as $\#\Gamma_\Sigma$.  
The Flexibility Theorem of Giroux \cite{Giroux91} states that 
$\Gamma_\Sigma$ (not the precise characteristic foliation) 
encodes all the contact-topological information in a small neighborhood of $\Sigma$.
The complement of the dividing set is the union of two subsets $\Sigma\setminus \Gamma_\Sigma=
\Sigma_+-\Sigma_-$.  Here $\Sigma_+$ is the subsurface where the orientation of $\Sigma$ and the normal orientation
of $\xi$ coincide, and $\Sigma_-$ is the subsurface where they are opposite. Therefore, we can refer to 
{\it positive} and {\it negative} components of $\Sigma\setminus \Gamma_\Sigma$.

\subsection{Symplectic fillings}\label{fillings}
The easiest way to prove a contact structure is tight is to show it `bounds' a symplectic 4-manifold.
There are several notions of `symplectic filling', and we assemble the various notions here for the
convenience of the reader.    (For more details, see the survey paper \cite{Etnyre98}.)

A symplectic manifold $(X,\omega)$ is said to have
{\em $\omega$-convex} boundary if there is a vector field $v$ defined in the neighborhood of $\partial X$ that
points transversely out of $X$ and for which $\mathcal{L}_v \omega=\omega,$ where $\mathcal{L}$ denotes the
Lie derivative. One may easily check that $\alpha=(\iota_v\omega)\vert_{\partial X}$ is a contact form on $\partial X.$
A symplectic manifold $(X, \omega)$ is said to have {\em weakly convex} boundary if $\partial X$
admits a contact structure $\xi$ such that $\omega\vert_\xi >0$  (and the orientations
induced on $\partial X$ by $X$ and $\xi$ agree).
A contact structure $\xi$ on $M$ is:

\be
\item {\it Holomorphically fillable} if $(M,\xi)$ is the $\omega$-convex boundary of some Stein manifold
$(X,\omega)$.
\item {\it Strongly symplectically fillable} if $(M,\xi)$ is the $\omega$-convex boundary of some
symplectic manifold $(X,\omega)$.
\item {\it Weakly symplectically fillable} if $(M,\xi)$ is the weakly convex boundary of some
symplectic manifold $(X,\omega)$.
\item {\it Weakly symplectically semi-fillable} if $(M,\xi)$ is one component of the weakly convex
boundary of some symplectic manifold $(X,\omega)$.
\ee

\begin{thm} [Gromov-Eliashberg] Let $(M,\xi)$ be a contact 3-manifold which satisfies any of the above conditions
for fillability.  Then $\xi$ is tight.
\end{thm}

The following diagram indicates the
hierarchy of contact structures.

$$
\begin{array}{|ccc|}
\hline
& &\\
\mbox{Tight} & &   \\
\cup \not | & &\\
\mbox{Weakly symplectically semi-fillable} & \supsetneqq & \mbox{Strongly symplectically
semi-fillable}\\ \cup  & & \cup \\
\mbox{Weakly symplectically fillable} & \supsetneqq  & \mbox{Strongly symplectically
fillable}\\ & & \cup\\
& & \mbox{Holomorphically fillable}\\
& &\\
\hline
\end{array}
$$

\s

The proper inclusion of the set of weakly symplectically semi-fillable contact structures into the set
of tight contact structures is the content of Theorem \ref{main}.  The proper inclusion of the set of
strongly fillable structures into the set of weakly fillable structures is already seen on $T^3$ (due to
Eliashberg \cite{E96}). This result was recently extended by Fan Ding to $T^2$-bundles over $S^1$.  For all
other inclusions in the diagram it is not known whether the inclusions are strict.

We briefly discuss when the various notions of fillability become the same.  If
$H^2(M;\Q)=0$, a weak symplectic filling can be modified into a strong symplectic filling
\cite{OO}.  
Using work of Kronheimer and Mrowka \cite{KM} and Seiberg-Witten theory, Lisca \cite{L} showed that 
if $M$ has a positive scalar curvature metric, then a semi-filling is
automatically a one-component filling 
(there is also a related, 
but weaker, result in Ohta-Ono \cite{OO}).


Lisca \cite{L} went further to show (among other things) that:

\begin{thm}[Lisca \cite{L}]\label{lisca}
Let $M$ be a Seifert fibered space over $S^2$ with Seifert invariants $(-{1\over 2},{1\over 3},
{1\over 4})$ or $(-{1\over 2}, {1\over 3}, {1\over 3})$.  The manifold $M$
does not carry a weakly symplectically semi-fillable contact structure.
\end{thm}

\begin{figure}[ht]
	{\epsfxsize=4.5in\centerline{\epsfbox{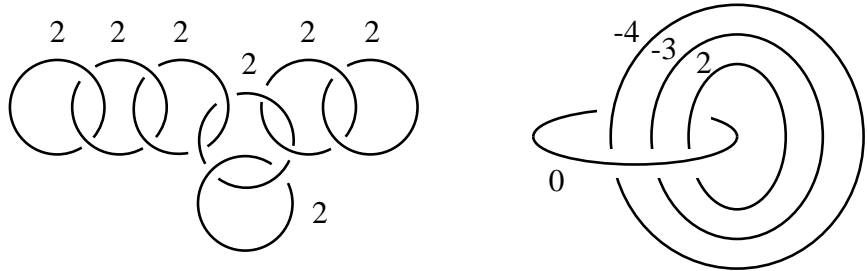}}}
	\caption{The plumbed disk bundles (left) and the Seifert fibered space $M$ (right).}
	\label{nofill}
\end{figure}

In Lisca's paper the Seifert fibered space with invariants $(-{1\over 2},{1\over 3}, {1\over 4})$ 
is described as the boundary $M$ of the 4-manifold
obtained by plumbing disk bundles over $S^2$ according to the positive $E_7$ diagram (left-hand
side of Figure \ref{nofill}).  It is an easy exercise in Kirby calculus \cite{GS} to show
that $M$ is orientation-preserving diffeomorphic to the manifold shown on the right-hand
side of Figure \ref{nofill}, which is  a presentation for a Seifert fibered
space.    Similarly, the Seifert fibered space with invariants $(-{1\over 2},{1\over 3},{1\over 3})$ 
corresponds to the positive $E_6$ diagram.   

\subsection{Legendrian surgeries}  \label{leg}
Let us now describe a contact surgery techinque called {\it Legendrian surgery}.  We first give a
description on the 3-manifold level.  Given a Legendrian knot $L$ in any contact 3-manifold
$(M,\xi)$, a {\em Legendrian surgery} on $L$ yields the contact manifold $(M',\xi')$,
where $M'$ is obtained from $M$ by $\tw(L)-1$ Dehn surgery on $L$ and
$\xi'$ is obtained from $\xi$ as follows: Let $N$ be a {\it standard convex neighborhood} of
$L$.  Choose a framing on $N$ so that $\tw(L)=0.$ This choice of framing allows us to
make an oriented identification $-\bdry(M\setminus N)\simeq \R^2/\Z^2$, where $(1,0)^T$ is
the meridian of $N$ and $(0,1)^T$ is the longitude of $N$ corresponding to the framing.
Now take an identical copy $N'$ of $N$ (with the same framing), and make an oriented
identification $\partial N'\simeq \R^2/\Z^2$, where $(1,0)^T$ is the meridian and
$(0,1)^T$ is the longitude. Then  let $M'=(M\setminus N)\cup_\psi N'$ where
$\psi:\partial N'\stackrel{\sim}{\rightarrow} -\partial (M\setminus N)$
is represented by the matrix $\left( \begin{array}{cc} 1 & 0 \\ -1 & 1\end{array} \right)\in
SL(2,\Z)$.  Since $\psi(\Gamma_{\partial N'})$ and $\Gamma_{\partial (M\setminus N)}$ are
isotopic, we may use Giroux's Flexibility Theorem to arrange the characteristic foliation on
$\bdry N'$ and isotop $\psi$ so $\psi_*(\xi|_{N'})=\xi|_{M\setminus N}$.  Hence we may glue the
contact structures on $N'$ and $M\setminus N$.

\begin{thm} \label{category-preserving} Legendrian surgery is category-preserving for each
category in the diagram of inclusions above, with the possible exception of the category of tight
contact structures. \end{thm}

Eliashberg \cite{Eliashberg90} proved that Legendrian surgery is category-preserving for
holomorphically fillable contact structures.  Weinstein \cite{Weinstein91} proved Theorem
\ref{category-preserving} for strongly symplectically fillable contact structures.

On the 4-manifold level, Legendrian surgery is described as follows:   Let $(X,\omega)$ be a
symplectic 4-manifold with $\omega$-convex boundary and $L$ a Legendrian knot in the induced
contact structure on $\partial X.$ If $X'$ is obtained from $X$ by adding a 2 handle to
$\partial X$ along $L$ with framing $\tw(L)-1,$ then $\omega$ extends to a symplectic form
$\omega'$ on $X'$ so that $\partial X'$ is $\omega'$-convex.
For many interesting applications of Theorem \ref{category-preserving}, we refer the reader to
Gompf \cite{Gompf}.

The case of Theorem \ref{category-preserving} known to a few experts but surprisingly absent in
the literature is for the category of weakly fillable contact structures.
The proof of Theorem \ref{category-preserving} in the strongly fillable case relies only on the
symplectic structure on $X$ in a neighborhood of $L$. Hence Theorem \ref{category-preserving}
in the weakly fillable case follows from:
\begin{lem}
	Let $(M,\xi)$ be a weakly symplectically fillable contact 3-manifold, $(X,\omega)$ one of its
	weak fillings, and $L$ a Legendrian knot in $(M,\xi).$ There is an arbitrarily small
	perturbation of $\xi$ in a neighborhood $N$ of $L$ so that $N$ is strongly convex. By this we
	mean there is a vector field 	$v$ defined on $X$ near $N$ so that $v$ points transversely out
	of $X,$  $\mathcal{L}_v\omega=\omega$ and $\xi\vert_N=\mbox{ker }(\iota_v\omega)\vert_N.$
\end{lem}
\begin{proof}
Let $(M',\xi')$ be any strongly fillable contact 3-manifold, $(X',\omega')$ one of its strong fillings and
$L'$ a Legendrian knot in $(M',\xi').$ It is not hard to find a neighborhood $N$ of $L$ in $M$ and $N'$ of $L'$ in $M'$
and a diffeomorphism $f:N\to N'$ such that $f(L)=L',$ $f^*\xi'=\xi$ along $L$, and 
$f^*(\omega'\vert_{N'})=\omega\vert_N.$ 
One may then use standard arguments
(see Exercise 3.35 in \cite{MS}) to extend $f$ to a symplectomorphism $(U,\omega)\stackrel{\sim}{\rightarrow} (U',\omega')$,
where $U\subset X$ is a neighborhood of $N$ and $U'
\subset X'$ is a neighborhood of $N'$. 
Finally note that the contact
planes $f^*\xi'$ and $\xi$ agree on $L$ and are close together near $L$.  
By a small perturbation of $\xi$ near $L$ (small enough to
keep $\xi$ contact so we may use Gray's Theorem), we may therefore assume $f^*\xi' = \xi$ near $L.$  Hence, if $v$ is the
expanding vector field for $\omega'$, then $f^{-1}_*(v)$ will be the desired vector field for $\omega$.
\end{proof}

We now comment on Theorem \ref{category-preserving} for the category of tight contact
structures.  In \cite{Honda4}, it was shown that there exists a tight contact structure on a
handlebody of genus 4 and a Legendrian surgery yielding an overtwisted contact structure.   It is
currently not known whether Legendrian surgery preserves tightness for {\it closed}
3-manifolds.

\section{The proof of the main result for the Seifert fibered space with invariants $(-{1\over 2}, {1\over 4}, 
{1\over 4})$}         \label{mainresult}

\subsection{Seifert fibered spaces}
Let $M$ be a Seifert fibered space over $S^2$ with three singular fibers $F_1, F_2, F_3$
and Seifert invariants $(\frac{\beta_1}{\alpha_1},\frac{\beta_2}{\alpha_2},
\frac{\beta_3}{\alpha_3})$.  We describe $M$ explicitly as follows: Let $V_i$, $i=1,2,3$, be a
tubular neighborhood of the singular fiber $F_i$.  We identify $V_i\simeq D^2\times S^1$
and $\bdry V_i\simeq \R^2/\Z^2$ by choosing $(1,0)^T$ as the meridional direction,
and $(0,1)^T$ as the longitudinal direction given by $\{pt\}\times S^1$.
We also identify $M\setminus(\cup_i
V_i)$ with $\Sigma_0\times S^1$, where $\Sigma_0$ is a sphere with three punctures, and further
identify $-\bdry (M\setminus V_i)=\R^2/\Z^2,$ by letting $(0,1)^T$
be the direction of an $S^1$-fiber, and $(1,0)^T$ be the direction given by $\bdry (M\setminus V_i)
\cap (\Sigma_0\times \{pt\})$.
With these identifications we may reconstruct $M$ from
$(\Sigma_0\times S^1)\cup(\cup_{i=1}^3 V_i)$
by gluing
$$A_i:\partial V_i \stackrel{\sim}{\rightarrow} -\partial (M\setminus V_i), \hspace{.15in}
 A_i=\left( \begin{array}{cc} \alpha_i & \gamma_i \\ -\beta_i & \delta_i\end{array}
\right)\in SL(2,\Z).$$   Note we have some freedom in choosing our matrices $A_i$ above. For
example, in choosing $A_i$ we can
alter  $\gamma_i$, $\delta_i$ by altering our choice of framing for $V_i$, which will
result in post-multiplying a given $A_i$ by $\left( \begin{array}{cc}
1 & m\\ 0 & 1\end{array} \right)$.

The twisting number of a Legendrian knot {\it isotopic} to
a regular ({\em i.e.}, non-singular) fiber of the Seifert fibration will be measured using the
framing from the product structure $\Sigma_0\times S^1$.  (Whenever we say {\it isotopy} we
will mean a smooth isotopy, as opposed to a {\it contact isotopy}.)  On the other hand, a
Legendrian knot isotopic to a singular fiber in a Seifert fibration will be measured with
respect to the framing on $V_i$ chosen in the description for $M.$



Let us now specialize to the case where $M$ is given by Seifert invariants $(-{1\over
2},{1\over 4},{1\over 4})$.  We make the following choices: $$A_1=
\left(\begin{array}{cc} 2 & -1 \\ 1 & 0\end{array}\right),
\hspace{.2in} A_2=A_3= \left( \begin{array}{cc}
4 & 1\\ -1 & 0\end{array} \right).$$
From now on $M$ will refer to this particular Seifert fibered space.

\subsection{Description as a torus bundle}
To define the contact structure $\xi$ in Theorem \ref{main} and prove tightness, we need a
description of $M$ as a torus bundle over $S^1$.
Recall a torus bundle over $S^1$ can be described as
$$T^2\times [0,1]/\sim,$$
where $(Ax,0)\sim(x,1)$ and $A$ is in $SL(2,\Z).$ The matrix $A$ is called the 
{\it monodromy} of the torus bundle.
\begin{lem}
	The manifold $M$ is a torus bundle over $S^1$ with monodromy $A=\left( 
	\begin{array}{cc} 0 & 1 \\ -1 & 0\end{array}
	\right).$
\end{lem}
\begin{proof}
The map $A$ has order four with two fixed points and two points of order two (which are interchanged
under $A$).
One may thus conclude
that the torus bundle is a Seifert fibered space over $S^2$ with Seifert invariants $(\pm\frac{1}{2},
\pm\frac{1}{4},\pm\frac{1}{4}).$ To determine the sign of invariants, let
$D\subset T^2$ be a small disk about one of the fixed points with $A(D)=D$, and consider
$S=D\times [0,1]/\sim$. If $x\in \bdry D$, then a regular fiber in the
Seifert fibered structure will be given by $(\{A^i(x)|i=0,1,2,3\}\times [0,1])/\sim$.
One may pick a product structure on $S$ so that a regular fiber will be a
$(-1,4)$-curve on $\partial S.$ From this we see that the gluing matrix $A_i$ associated to
this singular fiber is $\left( \begin{array}{cc} 4 & 1 \\ -1 & 0\end{array}
\right).$ Thus two of the Seifert invariants are $\frac{1}{4}.$ Similarly, one may check that
the third invariant is $-\frac{1}{2}.$
\end{proof}

\subsection{The tight contact structure}
We now define the tight contact structure $\xi$ on $M$, using the description of $M$ as a
torus bundle. For this we first describe a tight contact structure on $T^2\times [0,1]$ with
coordinates $(x,y,t)$. Start with a contact structure given by the 1-form  $\alpha=\sin({\pi
t\over 2})dx +\cos({\pi t\over 2})dy$. Perturb the boundary so that $T_i=T^2\times \{i\}$,
$i=0,1$, are convex with $\#\Gamma_{T_i}=2$, and slopes $s(\Gamma_{T_0})=0$,
$s(\Gamma_{T_1})=\infty$.   Now we identify $T_0$ and $T_1$ via $A$ to obtain a contact
structure on the quotient $M=(T^2\times [0,1])/\sim$.
This is possible since $A(\Gamma_{T_1})$ is isotopic to $\Gamma_{T_0}$ ---  we may apply
Giroux's Flexibility Theorem to ensure that the characteristic foliations agree.


\begin{prop}  \label{tightproof}
	The contact structure $\xi$ is a virtually overtwisted tight contact structure on $M.$
\end{prop}

In fact, $\xi$ is the {\it unique} virtually overtwisted tight structure on $M$ (see \cite{Honda2}). The uniqueness 
is not required in this paper.

\begin{proof}   The proof can be found in \cite{Honda2}, but we reproduce it here for
completeness.  We first show the existence of a double cover $\pi:M'\to M$ for which $\pi^*\xi$
is overtwisted.   If $M=(T^2\times[0,1])/\sim$, then
take two copies $C_1$ and $C_2$ of $(T^2\times [0,1],\xi)$ and map $T_1\subset
C_1$ to $T_0\subset C_2$ by $A$ and  $T_1\subset C_2$ to
$T_0\subset C_1$ by $A$.    In other words, $M'$ is the double cover of $M$ given
by $(T^2\times[0,2])/\sim$, where $(A^2x,0)\sim(x,2)$.   We may assume $\#\Gamma_{T_i}=2$,
$i=0,1,2$,  and $s(\Gamma_{T_0})=s(\Gamma_{T_2})=0$, $s(\Gamma_{T_1})=\infty$. If the
relative Euler class on $(T^2\times[0,1],\pi^*\xi)$ is  $e(T^2\times [0,1],
\pi^*\xi)=(0,1)-(1,0)=(-1,1)$, then $e(T^2\times[1,2],\pi^*\xi)=(1,1)$, and $e(T^2\times
[0,2],\pi^*\xi)=(0,2)$, which is not a possible relative Euler class for a tight contact
structure with the given boundary slopes. (See
\cite{Honda1} for a discussion of the relative Euler class.) Therefore $\pi^*\xi$ is overtwisted. 
We leave it as an exercise for the reader to explicitly find an overtwisted disk in $T^2\times [0,2]$.


We now prove the tightness of $\xi$. To construct $\xi$ on $M$, we started with a
tight contact structure $\xi|_{T^2\times [0,1]}$.  If there is an overtwisted disk
$D\subset M$, then necessarily $D\cap T_0\not=\emptyset$.   Below we show how to
inductively choose a different torus $T\subset M$ isotopic to $T_0$  that does
not intersect $D$, and for which $\xi|_{M\setminus T}$ is tight. This would contradict the
initial assumption of the existence of an overtwisted disk.

Without loss of generality we may assume that $D\pitchfork T_0$ and $D\cap T_0$ is a disjoint
collection of arcs and circles.
We first show how to eliminate an outermost arc of intersection.
Let $\alpha$ be an outermost arc of $D\cap T_0$ on $D$.  The arc $\alpha$ then separates off a
half-disk $D'$ from $D$ that does not intersect $T_0$ except along $\alpha$. If we cut open
$M$ to obtain $T^2\times [0,1]$, then $D'$ is attached to either $T_0$ or
$T_1.$ We assume the latter (the argument in the other case is similar). Let $N$ be a small
closed neighborhood of $T_1\cup D'$ in $T^2\times[0,1].$ The neighborhood $N$ is a
toric annulus which, by altering the product structure, we may call $T^2\times[\frac{1}{2},1].$
By taking $N$ to be small enough, we can ensure that $T_{1/2}=T^2\times\{\frac{1}{2}\}$ has
the same (isotopy class of) intersections with $D$ as $D\cap T_1$, except for the absence of
$\alpha.$
Moreover, we may assume that $T_{1/2}$ is convex. (For the moment we assume
that $\#\Gamma_{T_{1/2}}=2$. Below we show how to handle a convex
torus with more than two dividing curves.) Now glue $T_1$ to $T_0$ using
the map $A$ to obtain a toric annulus which we denote $T^2\times[-\frac{1}{2},\frac{1}{2}].$
Note that $(M,\xi)$ can be reconstructed by gluing the two boundary components of
$T^2\times[-\frac{1}{2},\frac{1}{2}]$, and $T_{1/2}\subset M$ is a
convex torus that has one fewer arc of intersection with $D$ than $D\cap T_0$.

Finally, we prove the contact structure $\xi|_{T^2\times[-1/2,1/2]}$ is tight.
By the classification of tight contact structures on toric annuli
\cite{Honda1, Giroux00}, any convex torus $T$ in $T^2\times[0,1]$ has slope $s(T)\leq 0.$
Let $s(T_{1/2})=-\frac{p}{q}<0$. If $s(T_{1/2})= 0$, then we
effectively did not modify the contact structure above.
$T^2\times[-\frac{1}{2},\frac{1}{2}]$ is a toric annulus with convex boundary having boundary
slopes $s(T_{-1/2})=\frac{q}{p}$ and $s(T_{1/2})=-\frac{p}{q}$ that was obtained by gluing
$T^2\times[-\frac{1}{2},0],$ a toric annulus with (universally) tight contact structure and
boundary slopes $\frac{q}{p}$ and 0, and $T^2\times[0,\frac{1}{2}],$ a toric annulus with
(universally) tight contact structure and boundary slopes 0 and $-\frac{p}{q}.$ Using the
classification of tight contact structures on toric annuli, we claim there is a tight contact
structure on $T^2\times[-\frac{1}{2},\frac{1}{2}]$ with boundary slopes $\frac{q}{p}$ and
$-\frac{p}{q}$ that may be split along the convex torus $T_0$ with $s(T_0)=0$, so
that the pieces are contactomorphic to the contact structures induced on
$T^2\times[-\frac{1}{2},0]$ and $T^2\times[0,\frac{1}{2}]$. This may easily be seen
using Theorem~4.24 in \cite{Honda1} which discusses the gluing of tight contact structures on
toric annuli. We have therefore eliminated one outermost arc of intersection from $D$ with a
convex torus in $M$ that cuts $M$ into a toric annulus with a tight contact structure. It is
not hard to continue this procedure to eliminate other outermost arcs of intersection.

We now eliminate circle components of $D\cap T_1.$ If $c$ is an
innermost circle of $D\cap T_1$ on $D$, then $c$ bounds disks $D'\subset D$ and $D''\subset
T_1$, and $D'\cup D''$ bounds a ball $B^3 \subset T^2\times[0,1]$. We may now apply the
same argument as above to a small neighborhood $N$ of $T_1\cup B$ to eliminate $c$ (and
possibly other curves) from $D\cap T_1$.  At each step of this reduction process, there exists
a convex torus $T$ which splits $M$ into a toric annulus with tight contact structure. It is
easy to see that by eliminating outermost arcs and innermost circles of $D\cap T$, we can
eventually make $D\cap T=\emptyset$, thus proving $\xi$ is tight.


Recall we must still discuss what happens if the new splitting torus $T_{1/2}$
constructed above has more than two dividing curves. In this case, Proposition~5.8 in
\cite{Honda1} implies there is a contact structure on a neighborhood
$T^2\times[\frac{1}{2}-\varepsilon,\frac{1}{2}+\varepsilon]$ of $T_{1/2}$ which is
$S^1$-invariant in the direction given by $s(T_{1/2})$ and contains a convex $T$ disjoint
from $T_{1/2}$, with $s(T)=s(T_{1/2})$  and $\#\Gamma_T=2$.
If we use $T$ as the new splitting torus, then $\xi|_{M\setminus T}$ will be tight (this is
just the argument above). But then $\xi|_{M\setminus T_{1/2}}$ will also be tight since there
exists a contact embedding
$$(M\setminus T_{1/2}, \xi|_{M\setminus T_{1/2}})\hookrightarrow (M\setminus T,
\xi|_{M\setminus T}).$$
\end{proof}

\subsection{Maximizing the twisting number of a regular fiber}

\begin{lem} \label{tw0}
There exists a Legendrian knot $F$ in $(M,\xi)$ isotopic to a regular fiber of the Seifert
fibered structure with $\tw(F)=0$.
\end{lem}

Lemma  \ref{tw0} follows immediately from the following lemma, together with
Proposition \ref{tightproof}.

\begin{lem}   \label{univtight}
	If $\xi'$ is a tight contact structure on $M$ and all Legendrian curves isotopic to
	a regular fiber have negative twisting number, then $\xi'$ is universally tight.
\end{lem}

\begin{proof}
To prove $\xi'$ is universally tight, we will first show that it can be made transverse to the
$S^1$-fibers of the Seifert fibration. We then show
that  this implies that $(M,\xi')$ is
covered by a tight contact structure on the 3-torus.  Since all of the tight contact structures
on the 3-torus are universally tight (see \cite{Kanda, Giroux00}),  $\xi'$ must therefore be
universally tight.

\s\n
{\bf Step 1.} (Normalization of contact structure $\xi'$.)
Let $F$ be a Legendrian curve isotopic to a regular fiber with $\tw(F)=n<0,$ which we take
to be maximal among Legendrian curves isotopic to a regular fiber. Let $L_i$, $i=1,2,3$,  be
Legendrian curves (simultaneously) isotopic to the singular fibers $F_i$ with $\tw(L_i)=n_i<0$,
and let $V_i$ be a standard convex neighborhood of $F_i$; assume also that $\tw(L_i)$ is
maximal among Legendrian curves isotopic to $F_i$ with negative twisting number.
After making the Legendrian ruling curves on $V_i$ vertical ({\it i.e.,} parallel to the
regular $S^1$-fibers), take a convex annulus $\mathcal{A}$ with Legendrian boundary for which
one component of $\bdry \mathcal{A}$ is a ruling curve on $V_2$ and the other component is a
ruling curve on $V_3$.   If not all dividing curves on $\mathcal{A}$ connect between $V_2$ and
$V_3$, then the Imbalance Principle (see \cite{Honda1}) gives rise to a bypass along a ruling
curve for, say, $V_2$.  Provided $\tw(L_2)<-1$, the Twist Number Lemma (see \cite{Honda1})
implies the existence of a Legendrian curve isotopic to $L_2$ with larger twisting number.
Therefore, we conclude that either $n_2=n_3<-1$ and $\mathcal{A}$ has no
$\bdry$-parallel dividing curves, or $n_2=n_3=-1$.

Assume $n_2=n_3\leq -1$  and $\mathcal{A}$ has no $\bdry$-parallel dividing curves.
If $N(\mathcal{A})$ is a convex neighborhood of $\mathcal{A}$, then $M'=V_2\cup V_3\cup
N(\mathcal{A})$ will have a piecewise smooth convex torus boundary. Rounding the corners in the
standard way (see \cite{Honda1}), $M'$ will be a convex torus with boundary slope
$${-n_2\over 4n_2+1}+{-n_2\over 4n_2+1}+{-1\over 4n_2+1}=-{2n_2+1\over 4n_2+1},$$
measured using the identification $\bdry (M\setminus V_1)\simeq \R^2/\Z^2$.  Now
$M''=\overline{M\setminus M'}$ is a solid torus with convex
boundary with slope $-{2n_2+1\over 4n_2+1}$, measured using $\bdry (M\setminus V_1)$, which is
equivalent to slope ${1\over 2n_2+1}$ measured using $\bdry V_1$.
This implies that the contact structure on $M''$ is the unique contact structure on the
standard neighborhood of the Legendrian knot $L_1$ with $n_1=2n_2+1$.

Now assume $n_2=n_3=-1$ and $\mathcal{A}$ has $\bdry$-parallel dividing curves.
Then we use the corresponding bypasses to thicken $V_2$, $V_3$ to $V_2'$, $V_3'$ so that the
boundary slopes (measured on $-\bdry (M\setminus V_i)$) are $-{1\over 2},-{1\over 2}$
or $-1,-1$ and the dividing curves on $\mathcal{A}$ between $V_2'$ and $V_3'$ do not have
$\bdry$-parallel curves.  The former case gives an overtwisted contact structure and the
latter yields a Legendrian curve isotopic to a regular fiber with $t=0$.

Summarizing, $\xi'$ has been normalized so that $V_2$ and $V_3$ are standard neighborhoods of
Legendrian curves with twisting number $n_2=n_3$, $\mathcal{A}$ has no $\bdry$-parallel
dividing curves, and $V_1=M\setminus (V_2\cup V_3\cup N(\mathcal{A}))$ is a standard
neighborhood of a Legendrian curve with twisting number $n_1=2n_2+1$.

\s\n
{\bf Step 2.} (Making $\xi'$ transverse to the fibers.)
Initially $K=\bdry V_2\cup \bdry V_3\cup \mathcal{A}$ has Legendrian rulings by vertical
curves. We perturb $K$ slightly so that the characteristic foliation becomes nonsingular
Morse-Smale, and $V_2\cap \mathcal{A}$ and $V_3\cap \mathcal{A}$ become transverse to $\xi'$.
Since $\bdry V_2$, $\bdry V_3$, and $\mathcal{A}$ are all convex in {\it standard form}, it is
possible to perturb $K$ along the Legendrian divides as in \cite{Honda2} to accomplish this.
Now, it is a question of isotoping $\xi'$ so that $\xi'$ is transverse to the fibers on
each $V_i$. Let us consider $V_2$, for example, and use the identification $\bdry V_2\simeq
\R^2/\Z^2$ to measure slope.  The regular fibers of the Seifert fibration have slope $-4$, and
the nonsingular Morse-Smale characteristic foliation has dividing curves of slope $-{1\over
n_2}$.  Since $-4<-{1\over n_2}$, it is clearly possible to extend $\xi'|_{\bdry V_2}$ so that
the contact structure is transverse to the Seifert fibers.  Moreover, this extension is contact
isotopic to $\xi'$ rel $\bdry V_2$.  In this way, we isotop $\xi'$ so that $\xi'$ is transverse
to the $S^1$-fibers of $M$.

\s\n
{\bf Step 3.} (Pulling back to $T^3$.)
Since the monodromy of $M$ as a torus bundle has order four, there exists a 4-fold
cover $\pi:T^3=T^2\times S^1\rightarrow M$ with the property that $\{pt\}\times S^1 \subset
T^3$ are lifts of fibers of $M$.
The pullback $\pi^*\xi'$ is therefore transverse to the fibers $\{pt\}\times S^1$ of $T^3$.

\s\n
{\bf Step 4.} (Universal tightness.)  
A tight contact structure on $T^3$ 
is universally tight by the classification results of Kanda \cite{Kanda} and Giroux \cite{Giroux00}.
Therefore, in order to show $\xi'$ is universally tight, it suffices to show that
$(T^3,\pi^*\xi')$ is tight.

We prove the well-known fact that a contact structure $\zeta$ on $T^2\times S^1$
which is transverse to the fibers must be tight. (The proof extends
easily to any circle bundle.) Let $\pi_1$, $\pi_2$ be the projections of $T^2\times D^2$ to
$T^2$ and $D^2$ respectively. If $\omega_1$, $\omega_2$ are area forms on $T^2$ and $D^2$, then
$\omega=\pi_1^*\omega_1+\pi_2^*\omega_2$ is a symplectic form
on $T^2\times D^2$ that restricts to a symplectic form on $\zeta$.  The contact structure $\zeta$ is
symplectically fillable and therefore tight.
\end{proof}

\subsection{Twisting number increase for singular fibers}

We need one more result before the proof of Theorem \ref{main}.
\begin{prop}\label{findl}
	There is a Legendrian knot $L$ in $(M,\xi)$ isotopic to one of the
	singular fibers $F_2$ or $F_3$ with $\tw(L)= 0$.
\end{prop}
\begin{proof}
Let $F$ be a Legendrian knot isotopic to a regular fiber with $\tw(F)=0$ as in Lemma \ref{tw0}.
Let $V_i'$, $i=1,2,3$, be disjoint solid tori isotopic to tubular neighborhoods of $F_i$, for
which $\bdry V_i'$ contains a contact-isotopic copy of $F$.  By perturbing $\bdry V_i'$ we may
assume $V_i'$ is convex with vertical ({\it i.e.,} parallel to the regular fibers) dividing curves, and,
furthermore, we may assume that $\#\Gamma_{\bdry V_i'}=2$, after possibly taking a smaller
solid torus.     In order to increase the twisting number of a Legendrian curve, we need to find a bypass. 
We will find a bypass along, say, $V_3'$  by patching together meridional disks of $V_1'$ and $V_2'$  to obtain
a punctured torus $T$ and showing the existence of a $\bdry$-parallel dividing curve on $T$.

\s\n
{\bf Step 1.}  (Normalizing $\xi$ on the complement.)  Let us first normalize the tight contact
structure on $\Sigma_0\times S^1 =M\setminus (V_1'\cup V_2'\cup V_3')$.

\begin{lem}\label{model}
The contact structure on $\Sigma_0\times S^1$ is contactomorphic to a
$[0,1]$-invariant tight contact structure on $T^2\times [0,1]$ with convex boundary,
$\#\Gamma_{T^2\times\{i\}}=2$, $i=1,2$, and  slopes $s(T^2\times \{i\})=\infty$, 
({\it i.e.,} the tight contact structure induced on $T^2\times[0,1]$, thought of
as a neighborhood of a convex torus in standard form) and a standard (open) neighborhood of a
vertical ({\it i.e.,} isotopic to $\{pt\}\times S^1\subset T^2$) Legendrian curve with 0 twisting
removed.
\end{lem}

\begin{proof}[Sketch of proof.]
This lemma is proved in \cite{EH} ({\em cf.\ }\cite{Honda2}), but we sketch the basic idea.
Since $\xi$ is tight on $M$, no dividing curve on $\Sigma_0$ can be $\bdry$-parallel.
This leaves us with two possibilities (depending on signs)  for the dividing curves on $\Sigma_0$, modulo
spiraling. Cutting $\Sigma_0\times S^1$ open along $\Sigma_0$, we obtain a tight contact
structure on a genus two handlebody with a fixed set of dividing curves. In addition the
dividing curves are such that one can use techniques in \cite{Kanda} to show there is a unique
tight contact structure on the handlebody. From this we can conclude that there is a unique
tight contact structure on $\Sigma_0\times S^1$ with the given dividing curve data. Now, since
the contact structure described in the lemma also has this dividing curve data, our contact
structure must be contactomorphic to it.
\end{proof}

\n
{\bf Step 2.} (Patching meridional disks.)
If we measure slopes of $\bdry V_i'$, $i=1,2,3$, using the identification $\bdry V_i\simeq
\R^2/\Z^2$, then the slopes are $2$, $-4$ and $-4,$ respectively.  After making the
ruling curves on $\partial V_i'$ meridional, a convex meridional disk $D_i$ for $V_i'$
will have, respectively, $\tb(\bdry D_i)=-2$, $-4$, $-4$, and also  2, 4 and 4 dividing curves.
We would like to patch copies of the meridional
disks together to create a convex surface and moreover relate information about the
dividing curves on this patched-together surface to the dividing curves on the meridional
disks.

We view the $T^2\times I$  (minus $S^1\times D^2$) from Lemma \ref{model} as the region between
$\partial V_1'$ and $\partial V_2'$ (minus $V'_3$).  Write 
$T_t=T^2\times \{t\}$, $t\in[0,1]$, as before.  Assume $T_0=\bdry V_2'$ and
$T_1=-\bdry V_1'$.   Since $\xi$ is $I$-invariant, we have (for example) a 1-parameter family
of positive regions $(T_t)_+=(T^2)_+\times\{t\}$.
We may then isotop $T_i$, $i=0,1$, away from $(T_i)_+$ ({\em i.e.,} on 
$(T_i)_-$) to arrange the slopes of the Legendrian ruling curves so that the meridional disk $D_i$
in $V_i'$ has Legendrian boundary. 
Now, take one copy of $D_2$ and
two copies $D_{11}$, $D_{12}$ of $D_1$, and arrange them so that $D_2\cap (T_0)_+=\delta\times \{0\}$
and $(D_{11}\cup D_{12})\cap (T_1)_+=\delta\times\{1\}$, where $\delta$ is a union of 
Legendrian arcs on $(T^2)_+$ with endpoints on opposite ends of $\bdry(T^2)_+$.  Let
$T=D_{11}\cup D_{12}\cup D_2\cup (\delta\times [0,1])$, which is a 
torus with an open disk removed.   See Figure~\ref{bypassonT}.
\begin{figure}[ht]
	{\epsfxsize=4in\centerline{\epsfbox{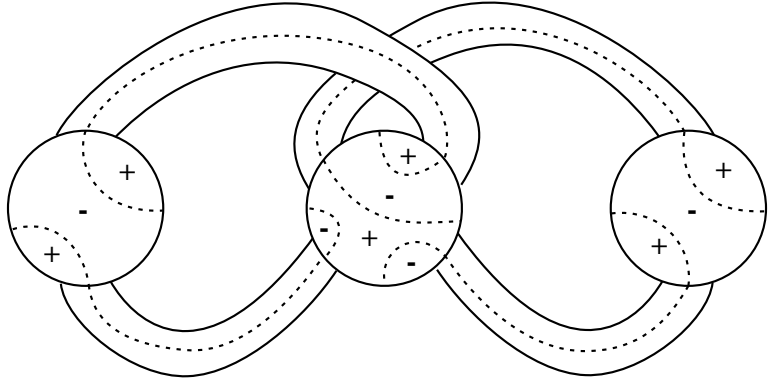}}}
	\caption{The punctured torus $T$ with dividing curves (dashed lines).} 
	\label{bypassonT} 
\end{figure}
After smoothing the corners using the `elliptic
monodromy lemma' or the `pivot lemma' of \cite{Etnyre, EF},  $T$ will have smooth Legendrian boundary. 
Since $\bdry T\subset \bdry ((T^2)_-\times [0,1])$,
we shall think of $T$ as having its boundary on $\partial V_3'.$

\s\n
{\bf Step 3.} (Combinatorics of $D_1$, $D_2$ and $D_3$.) 
Since $D_1$ has two dividing curves, it either has two positive regions and one negative region,
or one positive and two negative regions. We assume the former --- the argument for the latter 
is identical.

The rotation number $r(\bdry D_i)$, $i=2,3$, satisfies the formula    
$r(\bdry D_i)=\chi((D_i)_+)-\chi((D_i)_-)$ in \cite{Kanda98}.  Therefore, $r(\bdry D_i)$ can attain values $-3, -1, 1, 3$.


\s\n
{\bf Step 3A.} Assume that at least one of $D_2$ or $D_3$ (say $D_2$, after possible relabeling)
satisfies $r(\bdry D_i)>-3$. We first show that $D_2$, after possibly isotoping rel $\bdry D_2$,  
will have a positive $\bdry$-parallel region.  If $r(\bdry D_2)=3$ or $1$, there is no problem.  
If $r(\bdry D_2)=-1$, the dividing curves on $D_2$ may be either of the two types 
shown in Figure \ref{possible}.  
\begin{figure}[ht]
	{\epsfxsize=3.5in\centerline{\epsfbox{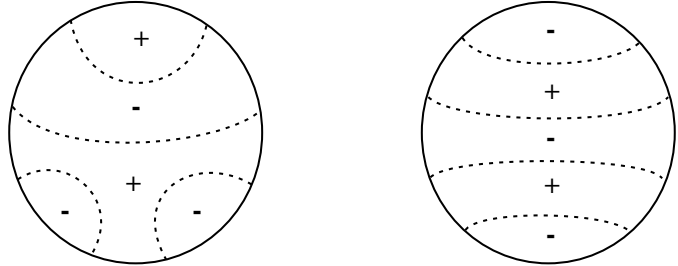}}}
	\caption{Possible dividing curves on $D_2.$}
	\label{possible}
\end{figure}
If we have a configuration shown on the right-hand side of 
Figure \ref{possible}, then we may isotop $D_2$ rel $\bdry D_2$ so that the dividing curves on 
$D_2$ are as shown on the left-hand side of Figure \ref{possible}.  This follows from  
the classification of tight contact structures
on solid tori in \cite{Honda1} or \cite{Giroux00}.

If $r(\bdry D_2)=-1$, then we take the dividing curves on $D_2$ to be as shown 
on the left-hand side of Figure~\ref{possible}. 
The dividing curves on $T$ will then  be as
shown in Figure~\ref{bypassonT}. 
Note we have a $\bdry$-parallel component and hence a bypass along $\partial T.$ 
The cases $r(\bdry D_2)=1,3$ are similar, by observing that any positive $\bdry$-parallel component of
$D_2$ must necessarily be connected to a positive $\bdry$-parallel component on
one of the copies of $D_1$, yielding a $\bdry$-parallel component on $T$.

The slope of $D_1$ on $\partial (M \setminus V'_1)$ is $-\frac{1}{2}$ and
the slope of $D_2$ on $\partial (M\setminus V'_2)$ is $\frac{1}{4}.$ This implies that
the slope of $T$ on $-\partial (M\setminus V'_3)$ is $-\frac{1}{4}.$ 
We  therefore have a bypass on $\partial V'_3$ attached along a ruling curve of slope 
$-\frac{1}{4}$ (as measured from
$\Sigma_0\times S^1$). Using this bypass, we may thicken $V'_3$ to $V''_3$ with standard convex boundary
having boundary slope $0.$ Thus, when measured from the product structure $D^2\times S^1$ on 
$V''_3$, the slope is
$\infty,$ showing that $V''_3$ is the standard neighborhood of a Legendrian curve $L$ isotopic to $F_3$ 
with twist number 0.

\s\n
{\bf Step 3B.} We are left with the case where $r(\bdry D_2)=r(\bdry D_3)=-3$. 
Now the dividing curves on
the punctured torus $T$ constructed from $D_2$ and two copies of $D_1$ will be as in Figure~\ref{tomany}. 
\begin{figure}[ht]
	{\epsfxsize=3in\centerline{\epsfbox{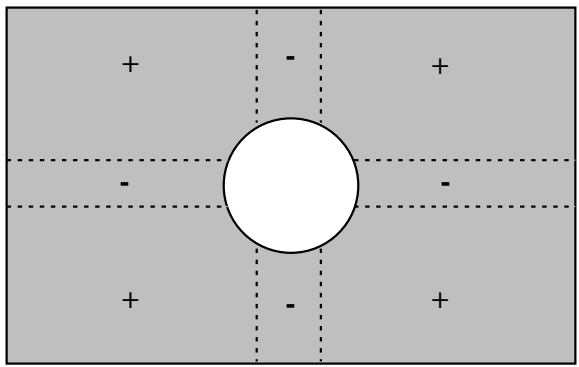}}}
	\caption{The punctured torus $T$ (shaded).}
	\label{tomany}
\end{figure}
Capping $T$
off with $D_3$, we obtain a closed contractible dividing curve on the torus $T\cup D_3$ which contradicts tightness.

This completes the proof of Proposition \ref{findl}. 
\end{proof}

\subsection{Proof of Theorem \ref{main} for the Seifert fibered space with invariants $(-{1\over 2},{1\over 4},{1\over 4})$. }
We are now ready to prove our main theorem.
\begin{proof}
From Proposition \ref{tightproof} we have a tight contact structure $\xi$ on $M.$ Now assume
that $\xi$ is weakly symplectically semi-fillable. From Proposition \ref{findl} we have a
Legendrian knot $L$ in $(M,\xi)$ that is isotopic to, say, $F_3$ with twist number 0.
We may assume that the neighborhood $V_3$ was chosen so that $L=F_3.$ As discussed above,
if we perform Legendrian
surgery on $L$, we remove a small neighborhood of $L$ (take this neighborhood to lie in
$V_3$) and re-glue it by $\left( \begin{array}{cc}
1 & 0\\ -1 & 1\end{array} \right).$ We easily see this has the same effect as changing
the $A_3$ to
$$\left( \begin{array}{cc}
3 & 1\\ -1 & 0\end{array} \right)=\left( \begin{array}{cc}
4 & 1\\ -1 & 0\end{array} \right)\left( \begin{array}{cc}
1 & 0\\ -1 & 1\end{array} \right).$$
Thus, after Legendrian surgery we obtain a weakly symplectically semi-fillable contact
structure on $M'$, the Seifert fibered space over $S^2$ with Seifert
invariants $(-\frac{1}{2},\frac{1}{3},\frac{1}{4}).$
This contradicts Lisca's Theorem (Theorem \ref{lisca}),
thus proving $\xi$ is not weakly symplectically semi-fillable.
\end{proof}

\subsection{Modifications for the Seifert fibered space with invariants $(-{2\over 3},{1\over 3},{1\over 3})$.}
The steps are almost identical for the Seifert fibered space $M_2$ over $S^1$ with 3 singular fibers
and invariants $(-{2\over 3},{1\over 3},{1\over 3})$. The manifold $M_2$ 
has a presentation as a torus bundle over $S^1$ with monodromy $A=\left( 
\begin{array}{cc}  
0 & 1\\ -1 & -1
\end{array}
\right)$.  By the classification in \cite{Honda2}, there exist two virtually overtwisted tight contact structures on 
$M_2$ which are nonisotopic but isomorphic. Let $\xi$ be one such contact structure.
As before, we first find a Legendrian knot  $F$ isotopic to 
a regular fiber with $t(F)=0$.  We then find a Legendrian curve $L$ isotopic to the 
$-{2\over 3}$-fiber $F_1$ with $t(L)$ large enough to perform a Legendrian surgery 
which modifies the Seifert invariants as follows:
$$\left(-{2\over 3}, {1\over 3},{1\over 3}\right)\rightsquigarrow \left(-{1\over 2}, {1\over 3}, {1\over 3}\right).$$
This again gives a contradiction of Theorem \ref{lisca}.  The existence of such an $L$ is proved by patching together
meridional disks as in Proposition \ref{findl} --- the only difference is that in one case we need to apply 
a `thinning before thickening'
argument that is used in \cite{EH}.




\section{Further Questions}

The obvious question raised in this paper is
\begin{question}
	Are the contact structures on the Seifert fibered spaces over $S^2$ with  Seifert
	invariants $(-{1\over 2},{1\over 3},{1\over 3})$ or $(-\frac{1}{2},\frac{1}{3},\frac{1}{4})$ constructed above tight?
\end{question}
We conjecture that these contact structures are tight. If the conjecture is true, we would have
an example of a manifold that supported a tight contact structure but no symplectically
fillable contact structures. Since this contact structure is constructed from a tight contact
structure by Legendrian surgery, we are led to ask the following:
\begin{question}
	Is Legendrian surgery category-preserving for tight structures on closed 3-manifolds?
\end{question}
If not, are there conditions which are sufficient to guarantee that Legendrian surgery on the
contact structure yields a tight contact structure?   For example,
\begin{question}
	Does Legendrian surgery on a universally tight contact structure produce a tight contact structure?
\end{question}
Recall our tight but not symplectically
fillable contact structure is virtually overtwisted. All other potential tight but not
symplectically fillable contact structures known to the authors are also virtually overtwisted.
(There are several candidates mentioned in \cite{Honda2}.)
So we ask
\begin{question}
	Are all universally tight contact structures symplectically semi-fillable?
\end{question}

\s\s\n
{\it Acknowledgements.}  We thank Yasha Eliashberg for informing us that Legendrian surgery 
preserves weakly fillable contact structures.


\begin{thebibliography}{99}

\bibitem{a}
B.\ Aebisher, et. al., \textit{Symplectic Geometry}, Progress in Math. \textbf{124},
	Birkh\"auser, Basel, Boston and Berlin, 1994.

\bibitem{Bennequin}
D.\ Bennequin, \textit{Entrelacements et \'equations de Pfaff},
	Ast\'erisque \textbf{107--108} (1983), 87--161.

\bibitem{Co97}
V.\ Colin, \textit{Chirurgies d'indice un et isotopies de sph\`eres dans les vari\'et\'es de contact
	tendues}, C.\ R.\ Acad.\ Sci.\ Paris Sér.\ I Math. \textbf{324} (1997), 659--663.

\bibitem{Co99}
V.\ Colin, {\it Recollement de vari\'et\'es de contact tendues},
	Bull. Soc. Math. France {\bf 127} (1999), 43--69.

\bibitem{Colin}
V.\ Colin, \textit{Chirurgie de Dehn admissible dans une vari\'et\'e de contact tendue}, preprint 2000.

\bibitem{Eliashberg89}
Y.\ Eliashberg, \textit{Classification of overtwisted contact structures on 3-manifolds},
	Invent. math. \textbf{98} (1989), 623--637.

\bibitem{Eliashberg90}
Y. Eliashberg, 
        \textit{Topological characterization of Stein manifolds of dimension $>$ 2},
        Int. J. of Math. \textbf{1} (1990), 29--46.

\bibitem{Eliashberg90a}
Y. Eliashberg, 
        \textit{Filling by holomorphic discs and its applications},
        Geometry of low-dimensional manifolds, Vol. II (Ed. Donaldson and Thomas),
        Cambridge, 1990.

\bibitem{Eliashberg92}
Y.\ Eliashberg, \textit{Contact 3-manifolds twenty years since J. Martinet's work},
	Ann. Inst. Fourier \textbf{42} (1992), 165--192.

\bibitem{E96}
Y.\ Eliashberg, {\it Unique holomorphically fillable contact structure on the $3$-torus},
	Internat. Math. Res. Notices {\bf 2} (1996),  77--82.


\bibitem{EF}
Y.\ Eliashberg and M. Fraser, \textit{Classification of topologically trivial 
	Legendrian knots}, in \textit{Geometry, topology, and
	dynamics} (Montreal, PQ, 1995), 17--51, CRM Proc.\ Lecture Notes
	\textbf{15}, Amer.\ Math.\ Soc., Providence, RI, 1998.

\bibitem{et}
Y.\ Eliashberg and W.\ Thurston, \textit{Confoliations},
	University Lecture Series \textbf{13}, Amer.\ Math.\ Soc., Providence, 1998. 

\bibitem{Etnyre98}  
J.\ Etnyre, \textit{Symplectic convexity in low-dimensional topology}, 
	Top. Appl. \textbf{88} (1998), 3--25.

\bibitem{Etnyre}
J.\ Etnyre, \textit{Contact structures on lens spaces}, 
     to appear in Commun. in Contemp. Math. 

\bibitem{EH} J.\ Etnyre and K.\ Honda, \textit{On the non-existence of tight contact structures},
	preprint 1999.

\bibitem{Giroux91}
E.\ Giroux, \textit{Convexit\'e en topologie de contact}, 
	Comment. Math. Helvetici \textbf{66} (1991), 637--677.

\bibitem{Giroux00}
E.\ Giroux, \textit{Structures de contact en dimension trois et bifurcations des feuilletages
	de surfaces}, Invent math. {\bf 141} (2000), 615--689.

\bibitem{Gompf} R.\ Gompf, \textit{Handlebody construction of Stein surfaces},
	Annals of Math.\ \textbf{148} (1998), 619--693.

\bibitem{GS} R.\ Gompf and A.\ Stipsicz, \textit{$4$-manifolds and Kirby calculus}, Graduate Studies in Mathematics,
	\textbf{20}, American Mathematical Society, Providence, RI, 1999.

\bibitem{Gromov}   M.\ Gromov, {\it Pseudo-holomorphic curves in symplectic manifolds},
	Invent. math. {\bf 82} (1985), 307--347.


\bibitem{Honda1} K.\ Honda, \textit{On the classification of tight contact
structures I}, Geom. Topol.  {\bf 4} (2000), 309--368.

\bibitem{Honda2} K.\ Honda, \textit{On the classification of tight contact
structures II}, to appear in J. Diff. Geom.

\bibitem{Honda4} K. Honda, \textit{Gluing tight contact structures}, preprint 2000.

\bibitem{Kanda} Y.\ Kanda, \textit{The classification of tight contact structures on
the 3-torus}, Comm. in Anal. and Geom. \textbf{5} (1997),  413--438.

\bibitem{Kanda98}
Y.\ Kanda, \textit{On the Thurston-Bennequin invariant of Legendrian knots and non exactness of
	Bennequin's inequality}, Invent.\ math.\ \textbf{133} (1998), 227--242.


\bibitem{KM} P.\ Kronheimer and T.\ Mrowka, \textit{Monopoles and contact structures}, 
	Invent.\ math.\ \textbf{130} (1997), 209--255.

\bibitem{L} P. Lisca, \textit{On symplectic fillings of $3$-manifolds}, Proceedings of 6th G\"okova Geometry-Topology
	Conference. Turkish J. Math. \textbf{23} (1999), 151--159.

\bibitem{LM} P. Lisca and G. Mati\'c, \textit{Stein 4-manifolds with boundary and contact structures},
	Top.\   Appl.\ \textbf{88} (1998), 55--66.

\bibitem{ML} S.\ Makar-Limanov, \textit{Morse surgeries of index 0 on tight manifolds}, preprint 1997.

\bibitem{MS}
D.\ McDuff and D.\ Salamon, \textit{Introduction to symplectic topology}, 
        Oxford University Press, 1995.

\bibitem{OO}
H.\ Ohta and K.\ Ono, \textit{Simple singularities and topology of symplectically filling 4-manifold},
	Comment.\ Math.\ Helv.\ \textbf{74} (1999), 575--590.

\bibitem{R}
L. Rudolph, \textit{An obstruction to sliceness via contact geometry and ``classical''
        gauge theory},
        Invent. math. \textbf{119} (1995), 155--163.

\bibitem{Weinstein91}
A. Weinstein, \textit{Contact surgery and symplectic handlebodies}, 
        Hokkiado Math. Journal \textbf{20} (1991), 241--251.


\end{thebibliography}
\end{document}